\documentclass [oneside]{amsart}
\usepackage {amscd}

\def\bc{\begin{center}}
\def\ec{\end{center}}
\def\px{\frac{\partial}{\partial x}}
\def\py{\frac{\partial}{\partial y}}
\def\pz{\frac{\partial}{\partial z}}

\def\G{\mathcal{G}}

\def\F{\mathcal{F}}

\newcommand{\C}{{\mathbb C}}
\newcommand{\Z}{{\mathbb Z}}
\newcommand{\Q}{{\mathbb Q}}
\newcommand{\N}{{\mathbb N}}

\newcommand{\CP}{{\mathbb {\mathbb{P}}^{2}}}

\newcommand{\cqd}{\ \hfill\rule[-1mm]{2mm}{3.2mm}}
\theoremstyle{change}
\newtheorem{theorem}{Theorem}
\newtheorem{conjecture}{Conjecture}
\newtheorem{proposition}{Proposition}
\newtheorem{lemma}{Lemma}
\newtheorem{definition}{Definition}

\newtheorem{corollary}{Corollary}
\newtheorem{example}{Example}
{ 
  \theoremstyle{margin}{
   }
}

\begin{document}

\title{On the Poincar\'e Problem for Foliations of General Type}
\author{Jorge Vit\'orio Pereira}

\maketitle
\bc
Instituto de Matem\'{a}tica Pura e Aplicada, IMPA, Estrada Dona Castorina, 110 \\
Jardim Bot\^{a}nico, 22460-320 - Rio de Janeiro, RJ, Brasil. 
email : jvp@impa.br  
\ec

\begin{abstract}
Let $\F$ be a holomorphic foliation of general type on $\CP$ which admits a rational first integral. We provide bounds for the degree of the first integral of $\F$ just in function of  the degree, the birational invariants of $\F$ and the geometric genus of a generic leaf. Similar bounds for invariant algebraic curves are also obtained and examples are given showing the necessity of the hypothesis.
\end{abstract}
\section{Introduction}

\rm
In \cite{Po}, Poincar\'e studied the following problem: "Is it possible to decide if an algebraic differential equation in two variables is algebraically integrable?". In the modern terminology the question above can rephrased as: "Is it possible to decide if a holomorphic foliation $\F$ on the complex projective plane $\CP$ admits a rational first integral ?".   Poincar\'e observed that in order to solve this problem is  sufficient to find a bound for the degree of the generic leaf of $\F$. For a modern version of some of Poincar\'e's results on this direction see \cite{Za}.

Since 1991 paper of Cerveau and Lins Neto \cite{CL},  this problem has received the attention of  many mathematicians. The interested reader should consult the papers by  Carnicer \cite{Ca}, Soares \cite{So1,So2}, Brunella--Mendes \cite{BM}, Esteves \cite{Es} and Zamora \cite{Za,ZZa} to have some idea of these recent developments.

Recently, see \cite{LN}, Lins Neto has constructed families of foliations with fixed degree and local analytic type of the singularities where foliations with rational first integral of arbitrarily large degree appear. This shows that it is impossible to bound the degree of the first integral in function of local information of the singularities. 

In this work we investigate the problem of bounding the degree of the first integral putting in evidence, for the first time in this problem, the Kodaira dimension of the foliation. This concept has been  introduced independently by L.G. Mendes \cite{LG} and  M. McQuillan \cite{Mc}. As in the birational theory of algebraic surface the Kodaira dimension measures the abundance of sections in the canonical bundle and its powers. For  algebraic surfaces the canonical bundle is nothing more than the bundle of holomorphic $2$--forms. In the foliated case it will be the bundle of $1$-forms defined over the leaves of the foliation. Formally, if we denote the canonical bundle of a holomorphic foliation $\F$ by $K_{\F}$ then the Kodaira dimension of $\F$ is defined as 
\[
  {\rm kod}(\F)= \limsup_{n \to \infty} \frac{ h^0(S, K_{\F}^{\otimes n}) }{\log n} \, ,
\]
and the sequence of integers $\lbrace P_n(\F)= h^0(S, K_{\F}^{\otimes n}) \rbrace_{n \in \N} $ is the {\it plurigenera} of $\F$.   

The possible values for the Kodaira dimension of a holomorphic foliation are $-\infty,0,1$ or $2$.  When ${\rm kod}(\F)=0$, $\F$ is, up to bimeromorphic maps and ramified coverings, generated by a global holomorphic vector field. If ${\rm kod}(\F) =1$ then $\F$ is a Riccati foliation or a turbulent foliation or some particular fibration. In the case $\F$ has  negative Kodaira dimension McQuillan has conjectured that $\F$ is a  rational fibration  or a  Hilbert modular foliation. Finally, when $\F$ is a foliation with ${\rm kod}(\F) =2$  we say that $\F$ is a foliation of {\it general type}. For more details see \cite{Mc}, \cite{LG} and \cite{Br}.

Our objective is to investigate the Poincar\'e problem for the foliations of general type. In this direction, our first result is:

\begin{theorem}\label{principal}
Let $\F$ be a holomorphic foliation of general type on $\CP$. Suppose that $\F$ admits a meromorphic first integral.  Then there exists a bound on the degree of the first integral depending only on the degree and the plurigenera of $\F$ and on the geometric genus of the generic leaf.
\end{theorem}

After some extra work we are able to extend the previous result to bound the degree of invariant curves of foliations which do not necessarily admit a rational first integral. More precisely we prove:

\begin{theorem}\label{segundo}
Let $\F$ be a holomorphic foliation of general type on $\CP$. Suppose that $\F$ admits an invariant algebraic curve $C$.  Then there exists a bound for the degree of $C$ depending only on the degree and the plurigenera of $\F$ and on the geometric genus of $C$.
\end{theorem}
  
Next we discuss some examples showing that the geometric genus has to appear as a parameter  for the bound, and that for foliations of Kodaira dimension distinct from two it is impossible to  bound the degree of invariant curves just in function of the degree of the foliation, even if the genus of the curve is fixed.

\section{The plurigenera of holomorphic foliations on surfaces}\label{revisao}
\rm
In this section we recall some basic definitions used in this work. Central to our exposition are the concepts of plurigenera and Kodaira dimension for holomorphic foliations. These concepts   were introduced in this  context  independently  by L.\ G.\ Mendes and M.\ McQuillan. For more information on the subject see \cite{Br}, \cite{LG} and \cite{Mc}.

A {\it holomorphic foliation} $\F$ on a compact complex surface $S$ is given by an open covering $\lbrace U_i \rbrace$ and holomorphic vector fields $X_i$ over each $U_i$ such that whenever the intersection of $U_i$ and $U_j$ is non--empty there exists an invertible holomorphic function $g_{ij}$ satisfying $X_i = g_{ij}X_j$. The collection $\lbrace g_{ij} \rbrace$ defines a holomorphic line--bundle $T{\F}$, called the {\it tangent bundle} of $\F$. The dual of $T{\F}$ is the {\it cotangent bundle} $T^*{\F}$, also called the {\it canonical bundle} $K_{\F}$.

\begin{definition}\rm
Let $\F$ be a foliation on the complex surface $S$. The {\it plurigenera} of $\F$ is defined as 
\[
  P_m(\F) = h^0(S, K_{\F}^{\otimes m}) \, \, \, \, \, {\rm for} \, \, \,  m \in \N^{+}.
\] 
\end{definition}

Recall that a {\it reduced foliation} $\F$ is a foliation such that every singularity $p$ is reduced in Seidenberg's sense, i.e., for every vector field $X$ generating $\F$ in a neighboorhoud of a singular point $p$, the eigenvalues of the linear part of $X$ are not both zero and their quotient, when defined, is not a positive rational number. 

\begin{definition}\rm
Let $\F$ be a  foliation on the complex surface $S$, and $\G$ any reduced foliation bimeromorphically equivalent to $\F$. The {\it Kodaira dimension} of $\F$ is given by
\[
  {\rm kod}(\F)= \limsup_{n \to \infty} \frac{ \log P_n(\G) }{\log n} \, .
\]
\end{definition}

When the foliation has Kodaira dimension $2$ we say that the foliation is of {\it general type}.

It follows from the next proposition that the Kodaira dimension is well defined and  is a bimeromorphic invariant of $\F$, for a proof see \cite{LG}.

\begin{proposition}\label{bi}
Let $\F$ be a reduced foliation on the algebraic surface $S_1$ and $\G$ be a reduced foliation on the algebraic surface $S_2$. If there exists a birational map $S_1 \to S_2$ sending $\F$ to $\G$, then $P_m(\F) = P_m(\G)$, for every nonnegative integer $m$. Consequently ${\rm kod}(\F) ={\rm kod}(\G)$.
\end{proposition}

\section{Bounds for the degree of a first integral}\label{enfim}

\noindent{\it \bf Proof of Theorem \ref{principal}}: Let $\sigma : S \to \CP$ be the minimal resolution of $\F$. Denote by $\G$ the reduced foliation $\sigma^*(\F)$. Since $\F$ admits a meromorphic first integral and $\G$ is free of dicritical singularities we have that $\G$ is a fibration. 

Let $C$ be a generic fiber of $\G$. Since $\F$ is of general type we can suppose that the genus $g$ of $C$ is at least $2$. Riemman-Roch Theorem implies that  
\[
  h^0(C, {\Omega^{1}_{C}}^{\otimes k})= k(2g -2) - g + 1 \, ,
\]
if $k$ is at least $2$. Take $n_0$ to be first non-negative integer that satisfies $P_{n_0}(\G) > n_0(2g -2) -g +1$. From the choice of $n_0$, the restriction map 
\[
        \phi : H^0(S,{{K_{\G}}}^{\otimes n_0}) \to H^0(C,{\Omega^{1}_{C}}^{\otimes n_0}) \, ,
\]
has  non empty kernel. In other words there exists a global holomorphic section of ${K_{\F}}^{\otimes n_0}$ that vanishes identically on $C$. 

Since $\sigma$ is a morphism, for every holomorphic section $s$ of ${{K_{\G}}}^{\otimes i}$ we have that $\sigma_*s$ is a holomorphic section of ${T\F ^*}^{\otimes i}$. Hence  
\[   
   \sigma_*(H^0(S,K_{\G}^{\otimes i})) \subset H^0(\CP,{T\F ^*}^{\otimes i}) \cong                H^0(\CP,\mathcal{O}_{\CP}(i \cdot {\rm d}(\F) -i)) \, ,
\]
one can conclude that the generic leaf of $\F$ has degree at most $n_0 \cdot ({\rm d}(\F) - 1) $.
\cqd

The bound obtained in the theorem above seems to depend on an infinite numbers of invariants of $\F$, $P_m(\F)$ for every positive integer $m$. Now we introduce a new bimeromorphic invariant, the {\it height} of a holomorphic foliation, which can be used to substitute the plurigenera as a parameter for the bound.

\begin{definition} \rm
Let $\F$ be a holomorphic foliation on an compact complex surface of non-negative Kodaira dimension $k$. As usual let $\G$ be any resolution of $\F$. We define the {\it height} of $\F$, $h(\F)$, to be the first positive integer $h$ such that  $K_{\G}^{\otimes h}$ has $k+1$ algebraically independent holomorphic sections.    
\end{definition}

\begin{lemma}\label{controla}
Let $\F$ be a reduced holomorphic foliation of general type on the compact complex surface $S$. If the height of $\F$ is $h$ then $P_{h \cdot n}(\F) \ge \binom{n + 2}{2}$.
\end{lemma}

\noindent{\it proof}: Let $V \subset H^0(S, K_{F}^{\otimes h})$ be a vector space generated by three algebraically independent global holomorphic sections of $K_{F}^{\otimes h}$. If we denote these sections by $s_0, s_1$ and $s_2$ and consider the morphism $\phi : S \to \CP$,
\[
   \phi(p) = (s_0(p) : s_1(p) : s_2(p) ) \, ,
\]
we obtain that $ V = \phi^* ( H^0(\CP,\mathcal{O}_{\CP}(1)) )$. Therefore we have that 
\[ 
     \phi^* (H^0(\CP, \mathcal{O}_{\CP}(n) )) \subset H^0(S, K_{F}^{\otimes h \cdot n }) \, ,
\]
 and the lemma follows. 
\cqd

\begin{corollary}
Let $\F$ be a foliation of general type on $\CP$. If $\F$ admits a meromorphic first integral  then there exists a bound for the degree of the first integral  depending only on the degree and the height of $\F$ and the geometric genus of the generic leaf.
\end{corollary}

\noindent{\it proof}: Follows easily from lemma \ref{controla} and the proof of theorem \ref{principal}. \cqd

A positive answer to the following conjecture would imply that the bound obtained in theorem \ref{principal} would depend just on the degree of $\F$ and the geometric genus of the generic leaf.

\begin{conjecture}
If $\F$ is a holomorphic foliation of general type on $\CP$ then there exists a bound for the height of $\F$ depending only on the degree of $\F$. 
\end{conjecture}

\section{Bounds for the degree of an invariant curve}

In order to extend Theorem \ref{principal} to bound degree of an invariant curve $C$ we must control the vanishing order of $\F$ along $C$. 

Let $p$ be  a reduced singularity of a holomorphic foliation $\F$ and $\Sigma$ be a local  smooth separatrix. If $\F$ in a neighboorhoud of $p$ is generated by a holomorphic vector field $X$ then the vanishing order of $\F$ along $\Sigma$ at $p$ is given by the Poincar\'e-Hopf index of $X_{|\Sigma}$ at $p$. We will use the notation $Z(\F,\Sigma,p)$. In more concrete terms, since $\Sigma$ is smooth at $p$ in a suitable coordinate system we can write
\[
   X_{|\Sigma} = \left( z^k + {\rm h.o.t.}  \right) \pz \, ,
\] 
and we set  $Z(\F,\Sigma,p)$  as $k$.

Note that this index is a particular case of the Gomez-Mont--Seade--Verjovski index which is defined for separatrices with arbitrary singularities.

\begin{lemma}\label{um}
Let $\F$ be a reduced holomorphic foliation on a surface $S$ and $C$ a smooth invariant curve. Then
\begin{enumerate}
\item  $\displaystyle{( K_{\F} )_{|C} = \Omega^{1}_C \otimes \mathcal{O}_C\left(\sum_{p \in C}{Z(\F,C,p) \cdot p}  \right)\,; }$
\item $h^0(C , ( K_{\F} )_{|C}^{\otimes m}) \le h^0(C,{\Omega^{1}_C}^{\otimes m}) + m \cdot Z(\F,C) $\, .
\end{enumerate}
\end{lemma}

\noindent{\it proof}: To prove item $1$ first consider an open covering $\lbrace U_i \rbrace$ of $S$ by Stein open sets. Over each $U_i$ of the covering take a holomorphic vector field $X_i$ generating $T\F$. Suppose that we have at most one singularity of $\F$ over each $U_i$. If the restriction of $X_i$ to $C$ does not have any singularity then we can interpret ${X_i}_{|C}$ as a local generator of $TC$, i.e., we have a canonical isomorphism between $( K_{\F} )_{|C\cap U_i}$ and $TC_{|U_i}$. When there exists a $p \in C$ such that $X_i(p) = 0$ then the restriction of $X_i$ to $C$ gives a section of $TC$ vanishing at $p$ with order $Z(\F,C,p)$,i.e., we have a canonical isomorphism between $( K_{\F} )_{|C\cap U_i}$ and $TC_{|U_i} \otimes \mathcal{O}_{U_i}(-Z(\F,C,p))$. 

Glueing the local canonical isomorphisms we obtain 
\[
  ( T{\F} )_{|C} = TC \otimes \mathcal{O}_C\left(\sum_{p \in   C} - {Z(\F,C,p) \cdot p}  \right)        \, ,
\]
and item $1$ follows by taking the dual.

Item $2$ follows  from item $1$ and the long exact sequence in cohomology associated to 
\[
  0 \to {\Omega^{1}_C}^{\otimes m} \to  {\Omega^{1}_C}^{\otimes m} \otimes \mathcal{O}_C (mD) \to \mathcal{O}_{mD} \to 0 \, ,
\]
where $D$ is the effective divisor given by 
\[
  D= \sum_{p \in C}{Z(\F,C,p) \cdot p }\, .
\]       \cqd

If $\F$ is a holomorphic foliation on $\CP$ we are going to say that $\G$ is the safe resolution of $\F$
if it is obtained by taking the minimal resolution of $\F$ and after that blowing-up each singularity 
once. We
do that in order  to guarantee that every irreducible  curve invariant by $\G$ is smooth.

\begin{lemma}\label{dois}
Let $\F$ be a holomorphic foliation on $\CP$ and $C$ an invariant algebraic curve. Let $\G$ be the safe resolution of $\F$ and $\overline{C}$ the strict transform of $C$. Then there exists a bound for $Z(\G,\overline{C})$ depending only on ${\rm d}(\F)$. 
\end{lemma}

\noindent{\it proof}: Let's say that a singularity is {\it quasi-reduced} if it has Milnor number $1$ or if it is a reduced saddle-node. First suppose that every singularity of $\F$ is quasi--reduced.

If $p$ is a reduced singularity then we have at most two branches of $C$ passing through $p$. After doing a blow-up at $p$ the contribution of the singularities infinitely near $p$ to $Z(\G,\overline{C})$ will be at most two. 

When $p$ is a dicritical singularity than we can have infinitely many branches of $C$ passing through $p$. Although after resolving $p$ it is not hard to see that at most two  branches of the strict transform of $C$ have singularities infinitely near $p$. Again the contribution of the singularities infinitely near $p$ to $Z(\G,\overline{C})$ is at most two.

If $p$ is a reduced saddle-node then the strong separatrix will contribute with one to $Z(\G,\overline{C})$ and the weak separatrix, if exists, with at most ${\rm d}(\F) + 1$.

Since the number of singularities of $\F$ is bounded by $d^2 +d + 1$, where $d= {\rm d}(\F)$,  the existence  of a bound for $Z(\G,\overline{C})$ in this particular case is proved.

To prove the general case one has just to observe two facts. 
The first fact is that there exists a positive integer $k$, depending on the degree of $\F$, such that with $k$ blow-ups we can always obtain a foliation  with all singularities quasi--reduced. 
The second fact is that a weak separatrix of a reduced saddle-node that appears after at most $k$ blow-ups will have bounded contribution to $Z(\G,\overline{C})$.
\cqd

\noindent{\it \bf Proof of Theorem \ref{segundo}}: The proof is completely similar to the proof of Theorem \ref{principal}. The only difference is that we have to use lemma \ref{um} and lemma \ref{dois} to guarantee the existence of the positive integer $n_0$ such that 
\[
  P_{n_0}(\F) > h^0(C,{{K_{\F}}_{|C}}^{\otimes n_0}) \, .
\] 
\cqd

Again we are able to substitute the plurigenera of $\F$ by the height of $\F$ and obtain the following corollary.

\begin{corollary}
Let $\F$ be a foliation of general type on $\CP$. If $\F$ admits an invariant algebraic curve $C$ then there exists a bound for the degree of $C$  depending only on the degree and the height of $\F$ and the geometric genus of $C$.
\end{corollary}

\section{Some Examples}\label{exemplos}

\begin{example}\label{negative}\rm
Let $(\F_{\alpha})_{\alpha \in \Q}$ the family of holomorphic foliations on $\CP$ defined, in an affine chart, by 
\[
   X_{\alpha} = x \px + \alpha y \py \, .
\]

Write $\alpha$ as $p/q$. If $\alpha$ is positive then we have a first integral of degree ${\rm max}(p,q)$, otherwise the first integral has degree $|p|+|q|$. Since the resolution of these foliations are rational fibrations, ${\rm kod}(\F_{\alpha}) = - \infty $.  
\qed
\end{example}

\begin{example}\label{zero}\rm
Let $(\F_{\alpha})_{\alpha \in \C}$ the family of holomorphic foliations on family of holomorphic foliations on $\CP$ defined by 
\[
   X_{\alpha} = (x^3 -1)(x - \alpha y^2) \px + (y^3 -1)(y - \alpha x^2) \py \, .
\]
These foliations, constructed by Lins Neto in \cite{LN}, are  foliations of degree $4$ on $\CP$ with  singularities with fixed local analytic type and with first integrals of arbitrarily large degree. Lins Neto also observed that whenever we have a first integral in this family then the generic leaf has geometric genus $1$. McQuillan, see \cite{Mc} or \cite{Br}, showed that Lins Neto examples up to a three-fold covering and birational transformations are linear vector fields on a torus and that the Kodaira dimension of each particular example is zero. 
\qed
\end{example}

\begin{example}\rm
The Gauss hypergeometric equation
\begin{equation}\label{chyper}
  z(1-z)w'' + \left( c - (a + b +1) z \right) w' - ab w = 0 \, ,
\end{equation}
whenever $c \notin \Z_{-}$, admits as general solution in the neighboorhoud of zero the function (see \cite{Hi})
\begin{equation}\label{csol}
\phi(z) = C_1 F(a,b,c;z) + C_2 z^{1-c}F(a-c+1,b-c+1,2-c; z)  \, ,
\end{equation}
where $C_1, C_2$ are arbitrary constants to be determined by the boundary conditions and
\[
  F(a,b,c;z) = 1 + \sum_{n=1}^{\infty} \frac{ \left( a \right)_n \left( b \right)_n  }{  \left( c \right)_n } z^n \, .
\]
Here $(p)_n$ stands for $p(p+1)(p+2)\cdots (p+n-1)$.

The classical change of variable $y(z) = - d \, {\rm log}\, w(z)$, see \cite{Hi} p. 104, associates a Riccati foliation to any second order differential equation. In this new coordinate the foliation induced by Gauss hypergeometric equation  can be written as 
\begin{equation}\label{eq1}
\omega = z(1-z) \, dy - \left( z(1-z)y^2  +  \left( c - (a + b +1) z \right) y + ab \right) \, dz \, .
\end{equation}

If $c \notin \Q$, $a = 1 -k$, $k \in \N$, and $b$ is arbitrary then the foliation induce by (\ref{eq1}) does not admit a rational first integral and  has an invariant rational curve of degree $k+1$ defined by the polynomial
\[
  y \cdot F(1-k,b,c;x) - F'(1-k,b,c; x) \, .
\]

Since these foliations have Kodaira dimension one, it is impossible to bound the degree of rational
curves (geometric genus $0$) for this class of foliations.
\qed
\end{example}

By pulling back the family presented in example \ref{zero} Lins Neto constructed holomorphic foliations  of general type on $\CP$,  again with  singularities with fixed local analytic type and with first integrals of arbitrarily large degree. 

\begin{example}\label{dois}\rm
Let $(G_{\alpha})_{\alpha \in \C}$ the family of holomorphic foliations on family of holomorphic foliations on $\CP$ defined by 
\[
   X_{\alpha} = (x^3 -1)(x - \alpha y^2) \px + (y^3 -1)(y - \alpha x^2) \py \, .
\]
It is shown in section $3.2$ of \cite{LN} that if we take the pull-back of $(G_{\alpha})_{\alpha \in \C}$ by the morphism $F : \CP \to \CP$ given in homogeneous coordinates by $F(X,Y,Z) = (X^r,Y^r,Z^r)$ , then the induced family of foliations, denoted by $(F_{\alpha})_{\alpha \in \C}$ will have the following properties.
\begin{itemize}
\item the foliations in the family have degree $3r+1$ ;
\item There is a finite set of parameters $A$ such that the restricted family $(F_{\alpha})_{\alpha \in \C \setminus A}$ has non degenerated singularities of fixed analytic type and between these singularities there are exactly $3r^2 + 6r + 3$ dicritical singularities;
\item There exists a countable and dense set of parameters $E \subset \C$, such that for any $\alpha \in E$ the foliation $\F_{\alpha}$ has a rational first integral of degree $d_{\alpha}$, satisfying the properties that for any $k>0$ the sets $\lbrace \alpha \in E | d_{\alpha} \le k \rbrace$ and
\[  \lbrace \alpha \in E | \text{ the geometric genus of the generic leaf is at most } k \rbrace 
\]
are finite.
\end{itemize} 

Suppose that $r$ is sufficiently large. If $\F_{\alpha}$ is a foliation in the family then  there exists a polynomial of degree $3r-1$ that  vanishes on all dicritical singularities of  
$\F_{\alpha}$. Hence the product of this polynomial with any linear homogeneous polynomial will lift to a section of the canonical sheaf of the resolved foliation. This is sufficient to assure that the foliation is of general type.

\qed
\end{example}

We finally remark that the bounds obtained for the degree of invariant algebraic curves on $\CP$  can be easily extended to any compact complex surface $S$ with Picard group isomorphic to $\Z$.

\section{Acknowledgements}
The author wants to thank C.\ Camacho, for many suggestions on the style of the exposition, and L.\ G.\ Mendes for many helpfull discussions and explanations about his work on the Kodaira dimension for holomorphic foliations. The author is supported by FAPERJ.



\begin{thebibliography}{99}

\bibitem {Br}
{\sc M.\ Brunella},
{\it Birational Geometry of Foliations},
First Latin American Congress of Mathematicians, IMPA, 2000.

\bibitem {BM}
{\sc M.\ Brunella and L.G. Mendes},
{\it Bounding the degree of solutions to Pfaff equations},
preprint 206, U. Bourgogne, 1999.


\bibitem{Ca}
{\sc M.\ Carnicer},
{\it The Poincar\'e problem in the nondicritical case},
Annals of Mathematics {\bf 140} (1994), 289--294.

\bibitem{CL}
{\sc D.\ Cerveau and A.\ Lins Neto},
{\it Holomorphic foliations in $\CP$ having an invariant algebraic curve},
Ann. Inst. Fourier {\bf 41} (1991), 883-903.

\bibitem {Es}
{\sc E.\ Esteves},
{\it The Castelnuovo-Mumford regularity of a variety left invariant by a vector field on projective space},
Preprint, IMPA, 2000.



\bibitem {Hi}
{\sc E.\ Hille},
{\it Ordinary differential equations in the complex domain},
John Wiley \& Sons, 1976.


\bibitem {LN}
{\sc A.\ Lins Neto },
{\it Some Examples for Poincar\'{e} and Painlev\'{e} Problems},
Preprint, IMPA, 2000.


\bibitem {LG}
{\sc L.\ G.\ Mendes},
{\it Kodaira dimension of holomorphic singular foliations},
Boletim da Sociedade Brasileira de Matem\'atica, {\bf 31}, 127--143, 2000

\bibitem {Mc}
{\sc M.\ McQuillan },
{\it Non--Commutative Mori Theory},
Preprint, IHES, 2000.

\bibitem {Po}
{\sc H. Poincar\'e},
{\it Sur l'integration alg\'{e}brique des \'{e}quations 
diff\'erentielles du premier ordre et du premier degr\'e} I and II, 
Rendiconti del Circolo Matematico di Palermo {\bf 5} (1891), 161--191; {\bf 11} (1897), 193--239.

\bibitem {Se}
{\sc J.\ B.\ Seaborn},
{\it Hypergeometric Functions and their applications},
Springer--Verlag, 1991.

\bibitem{So1}
{\sc M.\ Soares},
{\it The Poincar\'e problem for hypersurfaces invariant by one--dimensional foliations},
Invent. Math. {\bf 128} (1997), 495--500.

\bibitem{So2}
{\sc M.\ Soares},
{\it Projective varieties invariant by one--dimensional foliations},
to appear in Annals of Mathematics.


\bibitem{Za}
{\sc A.\ G.\ Zamora},
{\it Foliations in Algebraic Surfaces having a rational first integral},
Publicacions Matem\`{a}tiques {\bf 41} (1997), 357--373.

\bibitem{ZZa}
{\sc A.\ G.\ Zamora},
{\it Sheaves associated to holomorphic first integrals},
Ann. Inst. Fourier {\bf 50}, 3 (2000), 909--919 


\end{thebibliography}
\end{document}